\documentclass[12pt]{article}
\title{On the center of a compact group}
\author{Michael M\"uger\thanks{Supported by NWO through the `Pioneer' project no.\
616.062.384 of N.P. Landsman.} \\
Korteweg-de Vries Institute, Amsterdam, Netherlands \\ 
email: {\tt mmueger@science.uva.nl}}

\usepackage{latexsym, amssymb, amsmath, theorem}
\input diagrams         

\def\1#1{{\bf #1}}
\def\2#1{{\cal #1}}
\def\3#1{{\sl #1}}
\def\4#1{{\tt #1}}
\def\5#1{{\sf #1}}
\def\6#1{{\mathfrak #1}}
\def\7#1{{\mathbb #1}}

\newcommand{\ba}{\begin{array}}
\newcommand{\ea}{\end{array}}
\newcommand{\bea}{\begin{eqnarray}}
\newcommand{\eea}{\end{eqnarray}}
\newcommand{\bean}{\begin{eqnarray*}}
\newcommand{\eean}{\end{eqnarray*}}

\newcommand{\rarr}{\rightarrow}
\newcommand{\restr}{\upharpoonright}
\newcommand{\ol}{\overline}

\newcommand{\id}{\mathrm{id}}

\newcommand{\Hom}{\mathrm{Hom}}

\newcommand{\Obj}{\mathrm{Obj}\,}
\newcommand{\Rep}{\mathrm{Rep}}

\def\endexem{\hfill{$\Box$}\medskip}

\theoremstyle{change}
\theoremheaderfont{\scshape}
\newtheorem{defin}{Definition}[section]
\newtheorem{defprop}[defin]{Definition/Proposition}
\newtheorem{lemma}[defin]{Lemma}
\newtheorem{prop}[defin]{Proposition}
\newtheorem{theorem}[defin]{Theorem}
\newtheorem{coro}[defin]{Corollary}
\newtheorem{conj}[defin]{Conjecture}
\theorembodyfont{\rmfamily}
\newtheorem{rema}[defin]{Remark}
\newtheorem{noname}[defin]{}

\newcommand{\bdefin}{\begin{defin}}
\newcommand{\bdefprop}{\begin{defprop}}
\newcommand{\blemma}{\begin{lemma}}
\newcommand{\bprop}{\begin{prop}}
\newcommand{\btheor}{\begin{theorem}}
\newcommand{\bcoro}{\begin{coro}}
\newcommand{\edefin}{\end{defin}}
\newcommand{\edefprop}{\end{defprop}}
\newcommand{\elemma}{\end{lemma}}
\newcommand{\eprop}{\end{prop}}
\newcommand{\etheor}{\end{theorem}}
\newcommand{\ecoro}{\end{coro}}
\newcommand{\bconj}{\begin{conj}}
\newcommand{\econj}{\end{conj}}
\newcommand{\brem}{\begin{rema}}
\newcommand{\erem}{\endexem\end{rema}}
\newcommand{\bnix}{\begin{noname}}
\newcommand{\enix}{\end{noname}}

\newcommand{\prf}{\noindent{\it Proof. }}
\newcommand{\qed}{\ \hfill $\blacksquare$\medskip}

\begin{document}
\maketitle\noindent

\abstract{We prove a conjecture due to Baumg\"artel and Lled\'o \cite{BL} according to
which for every compact group $G$ one has $\widehat{Z(G)}\cong C(G)$, where the `chain
group' $C(G)$ is the free abelian group (written multiplicatively) generated by the set
$\widehat{G}$ of isomorphism classes of irreducible representations of $G$ modulo the
relations $[Z]=[X]\cdot[Y]$ whenever $Z$ is contained in $X\otimes Y$. Thus the
center $Z(G)$ depends only on the (ordered) representation ring of $G$. Furthermore, we
prove that every `t-map' $\varphi:\widehat{G}\rarr A$ into an abelian group, i.e.\
every map satisfying $\varphi(Z)=\varphi(X)\varphi(Y)$ whenever $X,Y,Z\in\widehat{G}$ and
$Z\prec X\otimes Y$, factors through the restriction map
$\widehat{G}\rarr\widehat{Z(G)}$. All these results generalize to pro-reductive groups
over algebraically closed fields of characteristic zero.}


\section{Introduction}
With every compact group $G$ one can associate two canonical compact abelian groups, to
wit the center $Z(G)$ and the abelianization $G_{ab}=G/\ol{[G,G]}$. Since every compact
group can be recovered from its (abstract) category of finite dimensional unitary
representations \cite{DR}, it is natural to ask whether the said abelian groups can be
recovered directly from $\Rep\,G$ without appealing to a reconstruction theorem \`a la
Tannaka-Krein-Doplicher-Roberts or Saavedra-Rivano-Deligne-Milne. Since $\Rep\,G$ is a 
discrete structure it is clear that one will rather recover the duals $\widehat{G_{ab}}$
and $\widehat{Z(G)}$. In the case of $\widehat{G_{ab}}$ it is well known how to proceed:
Writing $\2C=\Rep\,G$, let $\2C_1\subset\2C$ be the full subcategory of one dimensional
representations. Then the set of isomorphism classes of objects in $\2C_1$ is a (discrete)
abelian group, and it is easy to see that it is isomorphic to $\widehat{G_{ab}}$. It is
natural to ask whether also $\widehat{Z(G)}$ can be recovered directly from $\Rep\,G$. 

Motivated by certain operator algebraic considerations closely related to and inspired by
\cite{DR}, Baumg\"artel and Lled\'o \cite[Section 5]{BL} defined, for every compact group
$G$, a discrete abelian group $C(G)$ in terms of the representation category
$\Rep\,G$. They identified a surjective homomorphism $C(G)\rarr\widehat{Z(G)}$ and
conjectured the latter to be an isomorphism. They substantiated this conjecture by
explicit verification for several finite and compact Lie groups. (According to \cite{BL},
the case of $SU(N)$ was checked by C. Schweigert.) In this paper we prove
$\widehat{Z(G)}\cong C(G)$ for all compact groups,  exploiting a remark made in
\cite{eno}, and we derive two useful corollaries. Despite our general proof the examples
in \cite{BL} remain quite instructive.


\section{Definitions and Preparations}
Throughout the paper, $G$ denotes any compact group and $\widehat{G}$ the set of 
equivalence classes of irreducible representations. We allow ourselves the usual
harmless sloppiness of not always distinguishing between an irreducible representation $X$
and its equivalence class $[X]\in\widehat{G}$. (Thus `Let $X\in\widehat{G}$' means `Let 
$\2X\in\widehat{G}$ and let $X\in\Rep\,G$ be simple such that $[X]=\2X$'.)
While $\widehat{G}$ is a group iff $G$ is
abelian, there always is a notion of `homomorphism' into an abelian group:   

\bdefin \label{def-hom}
Let $G$ be a compact group and $A$ an abelian group. A map $\varphi: \widehat{G}\rarr A$
is called a t-map (tensor product compatible) if we have 
$\varphi(Z)=\varphi(X)\varphi(Y)$ whenever $X,Y,Z\in\widehat{G}$ and $Z\prec X\otimes Y$.
\edefin

\blemma 
If $\varphi:\widehat{G}\rarr A$ is a t-map then $\varphi(1)=1$, where the first $1$
denotes the trivial representation of $G$, and $\varphi(\ol{X})=\varphi(X)^{-1}$ for every
$X\in\widehat{G}$. 
\elemma

\prf We have $\varphi(1)=\varphi(1\otimes 1)=\varphi(1)\varphi(1)$, thus $\varphi(1)=1$. 
For any $X\in\widehat{G}$, we have $1\prec X\otimes \ol{X}$, thus
$1=\varphi(1)=\varphi(X)\varphi(\ol{X})$, which proves the second claim.  
\qed

The following proposition is essentially due to \cite{BL}:

\bprop \label{prop-univ}
 For every compact group $G$ there is a universal t-map $p_G: \widehat{G}\rarr C(G)$. 
(Thus for every t-map $\varphi: \widehat{G}\rarr A$ there is a unique homomorphism $\beta:
C(G)\rarr A$ of abelian groups such that 
\[ \begin{diagram}
   \widehat{G} & \rTo^{p_G} & C(G) \\ & \rdTo_{\varphi} & \dTo_{\beta} \\ && A
\end{diagram}\]
commutes.) Here the `chain group' $C(G)$ is the free abelian group (written
multiplicatively) generated by the set $\widehat{G}$ of isomorphism classes of irreducible
representations of $G$ modulo the relations $[Z]=[X]\cdot[Y]$ whenever $Z$ is contained in
$X\otimes Y$. The obvious map $p_G:\widehat{G}\rarr C(G)$ is a t-map.
\eprop

\prf We clearly must take $\beta$ to send the generator $[X]$ of $C(G)$ to $\varphi([X])$,
proving uniqueness. In view of the definition of a t-map this is is compatible with the
relations imposed on $C(G)$, whence existence of $\beta$.
\qed

\brem 1. The above elegant definition of $C(G)$ is due to J.\ Bernstein. In \cite{BL},
$C(G)$ was defined as $\widehat{G}/\sim$, where $\sim$ is the equivalence relation given
by $X\sim Y$ whenever there exist $n\in\7N$ and $Z_1,\ldots,Z_n\in\widehat{G}$ such that
both $X$ and $Y$ are contained in $Z_1\otimes\cdots\otimes Z_n$. Denoting the
$\sim$-equivalence class of $X$ is denoted by $\langle X\rangle$,\ $C(G)$ is an abelian
group w.r.t.\ the operations $\langle X\rangle\langle Y\rangle=\langle Z\rangle$, where
$Z$ is any irrep contained in $X\otimes Y$, and $\langle X\rangle^{-1}=\langle\ol{X}\rangle$. 
The easy verification of the equivalence of the two definitions is left to the reader.

2. A chain group $C(\2C)$, in general non-abelian, satisfying an analogous universal
property can be defined for any fusion category $\2C$, but we need only the case
$\2C=\Rep\,G$ and write $C(G)$ rather than $C(\Rep\,G)$. 
\erem

The following, proven in \cite{BL}, is the most interesting example of a t-map:

\bprop
The restriction of irreducible representations of $G$ to the center defines a surjective
t-map $r_G: \widehat{G}\rarr\widehat{Z(G)}$. Thus also the homomorphism of abelian groups 
$\alpha_G:C(G)\rarr\widehat{Z(G)}$ arising as above is surjective.
\eprop

\prf If $Z\in\widehat{G}$ and $g\in Z(G)$ then $\pi_Z(g)$ commutes with $\pi_Z(G)$, thus 
by Schur's lemma we have $\pi_Z(g)=\chi_Z(g)1_Z$, where $\chi_Z(g)\in U(1)$. Clearly,
$g\mapsto\chi_Z(g)$ is a homomorphism, thus $\chi_Z\in\widehat{Z(G)}$. This defines a map 
$r_G: \widehat{G}\rarr\widehat{Z(G)}$, which is easily seen to be a t-map. Since $Z(G)$ is
a closed subgroup of $G$, \cite[Theorem 27.46]{HR2} says that for every irreducible
representation (thus character) $\chi$ of $Z(G)$ there is a unitary representation $\pi$
of $G$ such that $\chi\prec\pi\restr Z(G)$. We thus have $r_G([\pi])=\chi$, thus $r_G$ is
surjective. Therefore also the map $\alpha_G: C(G)\rarr\widehat{Z(G)}$ arising from
Proposition \ref{prop-univ} is surjective.
\qed

For brevity we denote as fusion categories the semisimple $\7C$-linear tensor categories
with simple unit and two-sided duals, e.g.\ the $C^*$-tensor categories with conjugates,
direct sums and subobjects of \cite{DR}. (We do not assume finiteness.) All subcategories 
considered below are full, monoidal, replete (closed under isomorphisms) and closed under
direct sums, subobjects and duals, thus they are themselves fusion categories. 

\bdefin \label{def-C0}
Let $\2C$ be a fusion category. Then $\2C_0$ denotes the full tensor subcategory generated
by the simple objects $X$ for which the exists a simple object $Y\in\2C$ such that
$X\prec Y\otimes\ol{Y}$. 
\edefin

\brem The subcategory $\2C_0$ of a fusion category seems to have first been considered by
Etingof et al.\ \cite[Section 8.5]{eno}, where the following fact is remarked in
parentheses. The proof might be well known, but we are not aware of a suitable reference.
\erem

\bprop \label{prop-main}
Let $G$ be a compact group and $\2C=\Rep\,G$. Then the category $\2C_0$ coincides
with the full subcategory $\2C_Z\subset\2C$ consisting of those representations that are 
trivial when restricted to $Z(G)$. Thus $\2C_0\simeq\Rep(G/Z(G))$.
\eprop

\prf If $X,Y\in\widehat{G}$ are simple and $X\prec Y\otimes\ol{Y}$ then the restriction of
$X$ to $Z(G)$ is trivial, implying $\2C_0\subset\2C_Z$. As to the converse, let $g\in G$
be such that $g\in\ker\pi_X$ for all $X\in C_0$. This holds iff
$(\pi_Y\otimes\pi_{\ol{Y}})(g)=\11$ for all simple $Y\in\Rep\,G$. The latter means
\[ \pi_Y(g)\otimes\pi_Y(g^{-1})^t=\11, \]
which is true iff $\pi_Y(g)\in\7C\11_Y$. Now, if $g\in G$ is represented by scalars in 
all irreps $Y\in\widehat{G}$ then $g\in Z(G)$. (This follows from the fact that the
irreducible representations separate the elements of $G$.) In view of the Galois
correspondence of full monoidal subcategories $\2D\subset\Rep\,G$ and closed normal
subgroups $H\subset G$ given by 
\bean H_\2D &=& \{ g\in G \ | \ \pi_X(g)=\id\ \forall X\in\2D\}, \\
  \Obj\2D_H &=& \{ X\in \Rep\,G\ | \ \pi_X(g)=\id\ \forall g\in H\}. 
\eean
we have $H_{\2C_0}\subset Z(G)=H_{\2C_Z}$, thus $\2C_Z\subset\2C_0$ and therefore
$\2C_0=\2C_Z$. 
\qed

\blemma \label{lem-KG}
Let $G$ be compact and $\2C=\Rep\,G$. For a simple object $X\in\2C$ we have $p_G([X])=1$
iff $X\in\2C_0$.
\elemma

\prf If $Z$ and $X_i,Y_i,\ i=1,\ldots,n$ are simple with $X_i\prec Y_i\otimes\ol{Y}_i$ and
$Z\prec X_1\otimes\cdots\otimes X_n$ then 
$1,Z\prec Y_1\otimes\ol{Y}_1\otimes\cdots\otimes Y_n\otimes\ol{Y}_n$, thus $Z\sim 1$. This
implies that $p_G([X])=\langle X\rangle=1$ for every simple $X\in\2C_0$. Conversely, let
$X\in\2C$ be simple such that $p_G([X])=1$. This is equivalent to $X\sim 1$, thus there
are simple $Y_1,\ldots,Y_n$ such that $1,X\prec Y_1\otimes\cdots\otimes Y_n$. Then
$X\prec Y_1\otimes\cdots\otimes Y_n\otimes\ol{Y}_1\otimes\cdots\otimes\ol{Y}_n$, and
therefore $X\in\2C_0$.
\qed


\section{Results}
\btheor \label{theor-main}
The homomorphism $\alpha_G: C(G)\rarr\widehat{Z(G)}$ is an isomorphism for every compact
group $G$.
\etheor

\prf Since all maps in the diagram 
\[\begin{diagram}
  \widehat{G} & \rTo^{p_G} & C(G) \\ & \rdTo_{r_G} & \dTo_{\alpha_G} \\ && \widehat{Z(G)}
\end{diagram} \]
are surjective, $\alpha_G$ is an isomorphism iff $\ker p_G=\ker r_G$. By Lemma
\ref{lem-KG}, $[X]\in\ker p_G$ iff 
$X\in\2C_0$. On the other hand, $[X]\in\ker r_G$ iff $X\in\2C_Z$. By Proposition
\ref{prop-main} we have $\2C_0=\2C_Z$, thus we are done.  
\qed

$C(G)$ is defined in terms of the set $\widehat{G}$ and the multiplicities
$N_{ij}^k=\dim\Hom(\pi_k,\pi_i\otimes\pi_j),\ i,j,k\in\widehat{G}$ (the `fusion rules' in
physicist parlance). The same information is contained in the representation ring $R(G)$
provided we take its canonical $\7Z-$basis or its order structure \cite{han} into account.
We thus have the following

\bcoro \label{coro}
The center of a compact group $G$ depends only on the (ordered) representation ring
$R(G)$, not on the associativity constraint or the symmetry of the tensor category
$\Rep\,G$. (In general, both the latter are needed to determine $G$ up to isomorphism.)  
\ecoro 

\brem A considerably stronger result holds for {\it connected} compact groups: Every
isomorphism of the (ordered) representation rings of two such groups is induced by an
isomorphism of the groups, cf.\ \cite{han}. For non-connected groups this is wrong: The
finite groups $D_{8l}$ and $Q_{8l}$ are  non-isomorphic but have isomorphic representation
rings, cf.\ \cite{han}. Yet, as remarked in \cite[Section 5.1]{BL}, the centers are
isomorphic (to $\7Z/2\7Z$), as they must by Corollary \ref{coro}.
\erem

As an obvious consequence of Proposition \ref{prop-univ} and Theorem \ref{theor-main} we
have:

\bcoro \label{conj1}
Let $G$ be a compact group and $A$ an abelian group. Then every t-map 
$\varphi: \widehat{G}\rarr A$ factors through $\widehat{Z(G)}$, i.e.\ there is a
homomorphism $\beta: \widehat{Z(G)}\rarr A$ of abelian groups such that
\[ \begin{diagram}
   \widehat{G} & \rTo^{r_G} & \widehat{Z(G)} \\ & \rdTo_{\varphi} & \dTo_{\beta} \\ && A
\end{diagram}\]
commutes.
\ecoro

\brem This result should be considered as dual to the well known (and much easier) fact
that every homomorphism $G\rarr A$ from a group into an abelian group factors
through the quotient map $G\rarr G_{\mathrm{ab}}$.
\erem

\brem The results of this note were formulated for compact groups mainly because of the
author's taste and background. In view of \cite{del1} all results of this paper generalize
without change to pro-reductive algebraic groups over algebraically closed fields of
characteristic zero.
\erem

\noindent{\it Acknowledgments.} I am grateful to the DFG-Graduiertenkolleg `Hierarchie und
Symmetrie in mathematischen Modellen' for supporting a one week visit to the RWTH
Aachen. In particular I thank Fernando Lled\'o for the invitation, for drawing my
attention to \cite{BL}, many stimulating discussions, comments on the first version of
this paper and for the pizza.


\end{document}